\documentstyle[12pt]{article}
\def\IC{{\rm C}\kern-0.55em {^{_|}}\kern0.3em}

\def\dfrac{\displaystyle\frac}
\def\ba1{\begin{array}{l}}
\def\ea{\end{array}}

\def\x{\xi}
\def\y{\eta}
\def\z{\lambda}
\def\u{(u,\xi,\eta)}
\def\w{(u+v,\xi,\lambda)}
\def\v{(v,\eta,\lambda)}
\def\au#1{a_{#1}(u,\x,\y)}
\def\av#1{a_#1(v,\y,\z)}
\def\aw#1{a_#1(u+v,\x,\z)}
\def\MM#1{a_#1'(u,\x,\y)}
\def\bee{\begin{eqnarray}}
\def\ee{\end{eqnarray}}
\def\neqv{\equiv\!\!\!\!\!/\,}
\def\sn{sn(\lambda\,u +F(\x)-F(\y))} 
 
\def\cd{cd(\lambda\,u +F(\x)-F(\y))} 
\def\cd{cd(\lambda\,u +F(\x)-F(\y))} 
\def\tanh{tanh(\lambda\,u +F(\x)-F(\y))} 
\def\ll{(\lambda\,u +F(\x)-F(\y))}

\begin{document}
\begin{small}
\centerline{\bf\large  Classification of  Eight-vertex Solutions  of} 
\vskip6pt
\centerline{\bf\large  The Colored Yang-Baxter Equation }

\vskip10pt
\centerline{ Shi-kun WANG}
\vskip1cm

\begin{center}
\begin{minipage}{110mm}
{\bf Abstract:}\, In this paper 
all eight-vertex type solutions of the colored Yang-Baxter equation  
dependent on spectral as well as color parameter are given. 
It is proved that  they  are composed of three groups of basic solutions, 
three groups of their  degenerate forms and two groups of trivial solutions 
up to five 
solution transformations. Moreover, all non-trivial solutions can be classified 
into two types called Baxter type and Free-Fermion type. 
\end{minipage}
\end{center}
\vskip10pt

\subsubsection*{\large \S 0. \bf Introduction }

\quad

The Yang-Baxter or triangle equation which first appeared in Refs. 1-3 plays 
a prominent role in many branches of physics, for instance, 
in factorized $S$-matrices [4], exactly solvable models of statistical physics [5], complete integrable quantum and classical systems [6], quantum groups [7], conformal field theory and 
link invariants [8-11], to name just a few. In view of 
the importance of the Yang-Baxter equation, much attention has 
been directed to the search for solutions of the equation.

Colored Yang-Baxter equation dependent on spectral as well as color 
parameters is a generalization of usual 
Yang-Baxter equation. It has also attracted a lot of research interest 
(see Refs. 12-16) to find exact solutions for this type of Yang-Baxter 
equation. This is because the colored Yang-Baxter equation  concerns 
free Fermion model in magnetic field, multi-variable invariants of links  
and representations of quantum algebras and so on (see Refs. 17-21).  
 
The eight-vertex type solution of the colored Yang-Baxter 
equation has been investigated previously in Refs. 12, 17 and 20.
In Ref. 17 Fan and Wu first provided a single 
relation between the eight vertex weights in the general eight-vertex model 
by the pfaffian or  dimer method, the so-called free-fermion condition. 
Based on this work, V.V Bazhanov and Yu. G. Stroganov obtained a 
eight-vertex solution for the colored Yang-Baxter equation in Ref. 12,  
devoted to the eight-vertex free fermion model on a plane lattice. 
In Ref. 20 J. Murakami gave another eight-vertex solution in discussing 
multi-variable invariants of links. These are only two eight-vertex solutions   
for the colored Yang-Baxter equation we have known up to now.

The main theme of this paper is to give and classify  
all eight-vertex solutions  of the colored Yang-Baxter equation. The way 
to find the solution is from a computer algebra method given by Wu 
in Ref. 22. Moreover,  a theorem in the Ref. 22  can prove that all solutions 
can be obtained  by Wu's method.
The paper is organized as follows. In section one we  
will review the colored Yang-Baxter equation which in fact is a matrix 
equation. If the equation is expressed in components form it is composed of  
28 polynomial equations in the eight-vertex case.   
In this section we will first introduce the symmetries or solution 
transformation for this 
system of equations and some definitions of 
Hamiltonian coefficient, initial value condition, unitary condition and 
non-trivial gauge solution of the colored Yang-Baxter equation. Using the 
symmetries we can simply the system of equations as 12 polynomial equations.
In section two  we will apply the computer algebra method to the 12 
polynomial equations to get  
the  algebraic curves and differential equations satisfied by the
eight-vertex type solution of  the colored Yang-Baxter equation. In this  
section we will also give two relations satisfied by 
Hamiltonian coefficient which will play an important role in classification   
of eight-vertex solutions of the colored Yang-Baxter equation.
Based on the second  section, in the third
section we will  construct all non-trivial gauge eight-vertex type 
solutions of the colored Yang-Baxter equation and classify them 
into two types   called  Baxter and Free-Fermion type. The fourth section 
is devoted to general solutions and the relation between 
Hamiltonian coefficients and spin-chain Hamiltonian. In the third and fourth 
sections we will also show the two solutions appeared in Refs. 12 and 20 
are special cases of the general solutions obtained in this paper.

    In this paper symbolic computation will be applied 
to accomplish some tedious computations and the results obtained by 
computer calculations will be denoted by the symbol * . 

\vskip10pt
\subsubsection*{\large \S 1. \bf Colored Yang-Baxter equation, 
its symmetry and initial condition}

\quad

By  colored Yang-Baxter equation  we mean the following matrix equation
$$
\begin{array}{c}
\check{R}_{12}(u,\xi,\eta)\check{R}_{23}(u+v,\xi,\lambda) \check{R}_{12}(v,\eta,\lambda) =\check{R}_{23}(v,\eta,\lambda)\check{R}_{12}(u+v,\xi,\lambda) 
\check{R}_{23}(u,\xi,\eta),\\[4mm] 
\check{R}_{12}(u,\xi,\eta)=\check{R}(u,\xi,\eta)\otimes E,\qquad 
\check{R}_{23}(u,\xi,\eta)=E\otimes \check{R}(u,\xi,\eta),
\ea
\eqno (1.1)
$$
where $\check{R}\u$ is a matrix function of  $N^2$-dimension of $u, \x$ and
$\y$, $E$  is the unit matrix of order $N$ and  $\otimes$ means  tensor product 
of two matrices.  $u,v$ are called  spectral parameters and $\x$, $\y$ 
colored parameters. If the matrix is independent of colored parameters, 
then the colored Yang-Baxter equation (1.1) will become as the usual 
Yang-Baxter 
equation. If it is independent of spectral parameter then (1.1) will be 
reduced to the pure colored Yang-Baxter equation.
$$
\ba1
\check{R}_{12}(\xi,\eta)\check{R}_{23}(\xi,\lambda) \check{R}_{12}(\eta,\lambda) =\check{R}_{23}(\eta,\lambda)\check{R}_{12}(\xi,\lambda) 
\check{R}_{23}(\xi,\eta).
\ea\eqno(1.2)
$$

For the colored Yang-Baxter equation (1.1), the main interest in the paper is
to discuss the solutions with  the following form 
$$
\check{R}\u=\left(\begin{array}{cccc}
R^{11}_{11}\u  & 0             &0              &R^{11}_{22}\u\\
0              &R^{12}_{12}\u  &R^{12}_{21}\u  &0\\
0              &R^{21}_{12}\u  &R^{21}_{21}\u  &0\\
R^{22}_{11}\u &0              &0              &R^{22}_{22}\u
\end{array}\right).\eqno(1.3)
$$
The eight weight functions in (1.3) are denoted
by
$$
\begin{array}{ll}
\au1=R^{11}_{11}\u ,\quad \au5=R^{12}_{21}\u,\\[4mm]
\au2=R^{12}_{12}\u ,\quad \au6=R^{21}_{12}\u, \\[4mm]
\au3=R^{21}_{21}\u, \quad \au7=R^{11}_{22}\u,\\[4mm]
\au4=R^{22}_{22}\u, \quad \au8=R^{22}_{11}\u.
\ea
$$
We call the solution (1.3) eight-vertex type solution if  it satisfies
in addition $a_i(u,\xi,\eta)\neqv 0$ $(i=1,2,\cdots,8)$. Further, we only
consider the solutions which are meromorphic functions  of $u$ and $\x,\y$.
 For notational simplicity, throughout this paper let 
$$
u_i=a_i(u,\xi,\eta),\quad
v_i=a_i(v,\eta,\lambda),\quad
w_i=a_i(u+v,\xi,\lambda),\qquad
i=1,2,\cdots,8.
$$
For eight-vertex type solutions, the matrix equation (1.1) is equivalent to the following 28  equations:
\renewcommand{\theequation}{1.4\alph{equation}} \setcounter{equation}{0} \begin{eqnarray}
&\left .\ba1
u_{{7}}w_{{3}}v_{{8}}-u_{{8}}w_{{2}}v_{{7}}=0, \\[3mm] 
u_{{7}}w_{{8}}v_{{3}}-u_{{8}}w_{{7}}v_{{2}}=0, \\[3mm] 
u_{{2}}w_{{3}}v_{{2}}-u_{{3}}w_{{2}}v_{{3}}=0,  \\[3mm]
u_{{2}} w_{{8}}v_{{7}}-u_{{3}}w_{{7}}v_{{8}}=0,
\ea\right\}\qquad\\[4mm]
&\left .\ba1
u_{{1}}w_{{5}}v_{{2}}+u_{{7}}w_{{8}}v_{{6}}-v_{{2}}w_{{1}}u_{{5}}
-v_{{5}}w_{{2}}u_{{3}}=0, \\[3mm]
u_{{1}}w_{{1}}v_{{7}}+u_{{7}}w_{{3}}v_{{4}}-v_{{7}}w_{{5}}u_{{5}} 
-v_{{1}}w_{{7}}u_{{3}}=0, \\[3mm]
u_{{2}}w_{{6}}v_{{1}}+u_{{5}}w_{{7}}v_{{8}}-v_{{6}}w_{{1}}u_{{2}}
-v_{{3}}w_{{2}}u_{{6}}=0, \\[4mm]
u_{{1}}w_{{2}}v_{{1}}+u_{{7}}w_{{4}}v_{{8}}-v_{{2}}w_{{1}}u_{{2}}
-v_{{5}}w_{{2}}u_{{6}}=0, \\[3mm]
u_{{1}}w_{{7}}v_{{5}}+u_{{7}}w_{{6}}v_{{3}}-v_{{7}}w_{{5}}u_{{2}}
-v_{{1}}w_{{7}}u_{{6}}=0, \\[3mm]
u_{{1}}w_{{7}}v_{{2}}+u_{{7}}w_{{6}}v_{{6}}-v_{{1}}w_{{1}}u_{{7}}
-v_{{7}}w_{{2}}u_{{4}}=0, 
\ea\right\}\qquad\\[4mm]
&\left .\ba1
u_{{4}}w_{{6}}v_{{2}}+ u_{{7}}w_{{8}}v_{{5}}- v_{{2}}w_{{4}}u_{{6}}
-v_{{6}}w_{{2}}u_{{3}}=0, \\[3mm]
u_{{4}}w_{{4}}v_{{7}}+ u_{{7}}w_{{3}}v_{{1}}- v_{{7}}w_{{6}}u_{{6}} 
-v_{{4}}w_{{7}}u_{{3}}=0, \\[3mm]
u_{{2}}w_{{5}}v_{{4}}+ u_{{6}}w_{{7}}v_{{8}}- v_{{5}}w_{{4}}u_{{2}}
-v_{{3}}w_{{2}}u_{{5}}=0, \\[4mm]
u_{{4}}w_{{2}}v_{{4}}+ u_{{7}}w_{{1}}v_{{8}}- v_{{2}}w_{{4}}u_{{2}}
-v_{{6}}w_{{2}}u_{{5}}=0, \\[3mm]
u_{{4}}w_{{7}}v_{{6}}+ u_{{7}}w_{{5}}v_{{3}}- v_{{7}}w_{{6}}u_{{2}}
-v_{{4}}w_{{7}}u_{{5}}=0, \\[3mm]
u_{{4}}w_{{7}}v_{{2}}+ u_{{7}}w_{{5}}v_{{5}}- v_{{4}}w_{{4}}u_{{7}}
-v_{{7}}w_{{2}}u_{{1}}=0, 
\ea\right\}\qquad\\[4mm]
&\left .\ba1
u_{{1}}w_{{5}}v_{{3}}+ u_{{8}}w_{{7}}v_{{6}}- v_{{3}}w_{{1}}u_{{5}}
-v_{{5}}w_{{3}}u_{{2}}=0, \\[3mm]
u_{{1}}w_{{1}}v_{{8}}+ u_{{8}}w_{{2}}v_{{4}}- v_{{8}}w_{{5}}u_{{5}} 
-v_{{1}}w_{{8}}u_{{2}}=0, \\[3mm]
u_{{3}}w_{{6}}v_{{1}}+ u_{{5}}w_{{8}}v_{{7}}- v_{{6}}w_{{1}}u_{{3}}
-v_{{2}}w_{{3}}u_{{6}}=0, \\[4mm]
u_{{1}}w_{{3}}v_{{1}}+ u_{{8}}w_{{4}}v_{{7}}- v_{{3}}w_{{1}}u_{{3}}
-v_{{5}}w_{{3}}u_{{6}}=0, \\[3mm]
u_{{1}}w_{{8}}v_{{5}}+ u_{{8}}w_{{6}}v_{{2}}- v_{{8}}w_{{5}}u_{{3}}
-v_{{1}}w_{{8}}u_{{6}}=0, \\[3mm]
u_{{1}}w_{{8}}v_{{3}}+ u_{{8}}w_{{6}}v_{{6}}- v_{{1}}w_{{1}}u_{{8}}
-v_{{8}}w_{{3}}u_{{4}}=0, 
\ea\right\}\qquad\\[4mm]
&\left .\ba1
u_{{4}}w_{{6}}v_{{3}}+ u_{{8}}w_{{7}}v_{{5}}- v_{{3}}w_{{4}}u_{{6}}
-v_{{6}}w_{{3}}u_{{2}}=0, \\[3mm]
u_{{4}}w_{{4}}v_{{8}}+ u_{{8}}w_{{2}}v_{{1}}- v_{{8}}w_{{6}}u_{{6}} 
-v_{{4}}w_{{8}}u_{{2}}=0, \\[3mm]
u_{{3}}w_{{5}}v_{{4}}+ u_{{6}}w_{{8}}v_{{7}}- v_{{5}}w_{{4}}u_{{3}}
-v_{{2}}w_{{3}}u_{{5}}=0, \\[4mm]
u_{{4}}w_{{3}}v_{{4}}+ u_{{8}}w_{{1}}v_{{7}}- v_{{3}}w_{{4}}u_{{3}}
-v_{{6}}w_{{3}}u_{{5}}=0, \\[3mm]
u_{{4}}w_{{8}}v_{{6}}+ u_{{8}}w_{{5}}v_{{2}}- v_{{8}}w_{{6}}u_{{3}}
-v_{{4}}w_{{8}}u_{{5}}=0, \\[3mm]
u_{{4}}w_{{8}}v_{{3}}+ u_{{8}}w_{{5}}v_{{5}}- v_{{4}}w_{{4}}u_{{8}}
-v_{{8}}w_{{3}}u_{{1}}=0.
\ea\right\}
\end{eqnarray}

Assume $\check{R}(u,\xi,\eta)$ is a solution of (1.1). 
Having carefully studied the  system of equations (1.4), we find 
there are five symmetries for  eight-vertex type solutions of  the colored
Yang-Baxter equation (1.1).
\vskip5pt

 {\bf(A)\, Symmetry of interchanging indices.} 

The system of  equations (1.4) is invariant if we interchange
the two sub-indices $2$ and $3$ as well as the two sub-indices $7$ and $8$
or the sub-indices $1$ and $4$ as well as the two sub-indices $5$ and $6$. 
\vskip5pt

{\bf(B)\, The scaling symmetry.}

Multiplication of the solution $\check{R}\u$ by an arbitrary function
$g\u$ is still a solution of the colored Yang-Baxter equation (1.1).
\vskip5pt

{\bf (C)\, Symmetry of weight functions} 

If the weight functions $ a_2(u,\xi,\eta), a_3(u,\xi,\eta),a_7(u,\xi,\eta)$ 
and $a_8(u,\xi,\eta)$ are replaced by the new weight functions
$$ 
\ba1
\bar{a}_2(u,\xi,\eta)=\dfrac{N(\xi)}{N(\eta)}a_2(u,\xi,\eta), \quad
\bar{a}_3(u,\xi,\eta)=\dfrac{N(\eta)}{N(\xi)}a_3(u,\xi,\eta),\\[4mm] 
\bar{a}_7(u,\xi,\eta)=s\, N(\eta)N(\xi)a_7(u,\xi,\eta), \quad
\bar{a}_8(u,\xi,\eta)=\dfrac{1}{s\, N(\eta)N(\xi)}a_8(u,\xi,\eta) 
\ea 
$$
respectively or $\au5$ and $\au6$ by $-\au5$ and $-\au6$, where $N(\xi)$ is 
an arbitrary function of colored parameter
and $s$ is a complex constant, the new matrix $\check{R}(u,\xi,\eta)$ 
is still a solution of (1.1).
\vskip5pt
{\bf (D)\, Symmetry of spectral  parameter.}

If we take a new spectral parameter $\bar{u}=\mu u$ 
where $\mu$ is a complex constant independent of spectral and colored parameters
, the new matrix $\check{R}(u',\x,\y)$
is still a solution of (1.1).
\vskip5pt

{\bf (E)\, Symmetry of color parameters.}
If we take 
 new colored parameters $\zeta=f(\xi)$, $\theta=f(\eta)$, 
where  $f(\xi)$ is an arbitrary function
, then the new matrix $\check{R}(u,\zeta,\theta)$
is also a solution of (1.1).
\vskip5pt

The five symmetries (A)--(E) are called 
solution transformations A--E of eight-vertex type solutions
of  the colored Yang-Baxter equation (1.1) respectively.
\vskip5pt

Dividing  both sides of the third of   (1.4a) by
$\au2\aw2\av2$, we get 
$$
f\w=f\u f\v, 
\eqno (1.5)
$$
where $f\u=\au3/\au2$.  Putting $u=v=\y=0$ in (1.5) we have
$$
 f(0,\x,\z)=f(0,\x,0)f(0,0,\z).
$$
Substituting this formula into the one obtained by taking $u=v=\x=0$ in (1.5) 
we get
$$
f(0,0,\z)=f(0,0,\y)f(0,\y,\z)=f(0,0,\y)f(0,\y,0)f(0,0,\z).
$$
This means 
$$
f(0,\y,0)f(0,0,\y)=1.
$$
Otherwise, it is easy to show that $f\u = 0$, i. e. $\au3 = 0$.
Therefore 
$$
f(0,\x,\y)=\dfrac{M(\x)}{M(\y)},
\eqno (1.6)
$$ 
where $M(\x)=f(0,\x,0)$. On the other hand, if we differentiate both sides 
of (1.5) with respect to the spectral variable $v$ and then set 
$v=0,\,\z=\y$, then  
$$
f'(u,\x,\y)=f'(0,\y,\y)f(u,\x,\y)
\eqno(1.7a)
$$
holds, where the dot means  derivative with respect to $u$ and the 
simple formula 
$$
\dfrac{d\,G(u+v)}{d\,v}\big |_{v=0}=\dfrac{d\,G(u)}{d\,u},
$$
for any function $G(u)$, is used. Similarly, one also has 
$$
f'(v,\x,\z)=f'(0,\x,\x)f(v,\x,\z)
\eqno(1.7b)
$$
if we differentiate both sides of (1.5) with respect to $u$ and then set $u=0$
and $\y=\x$. The two formulas above imply $ f'(0,\x,\x)$ is a constant 
independent of  colored parameter $\x$.  Hence
$$
f\u=\dfrac{M(\x)}{M(\y)}\exp(\nu u), \eqno(1.8)
$$
where $\nu$ is a complex constant.

From the group of equations (1.4a) we have
\renewcommand{\theequation}{1.9\alph{equation}}
\setcounter{equation}{0}
\bee
&f\u h\v=h\w,\\
&f\u f\v=f\w,\\
&h\u =h\w f\v,\\
&h\u =f\w h\v,
\ee
where $ f\u=a_3(u,\x,\y)/a_2(u,\x,\y),\, h\u=a_8(u,\x,\y)/a_7(u,\x,\y)$. 
If we let $v=0$ and $\z=0$ in (1.9d) then 
$$
h\u=f(u,\x,0)h(0,\y,0).
$$
Substituting this into (1.9c) and using (1.9b), one obtains
$$
f(u,\x,0)h(0,\y,0)=f(u,\x,0)f(v,0,0)h(0,\z,0)f(v,\y,\z)
$$
or
$$
h(0,\y,0)=f(v,0,0)h(0,\z,0)f(v,\y,\z).
$$ 
This formula implies $\nu=0$ in (1.8) and then one can obtain
$$
f\u=\dfrac{M(\x)}{M(\y)},\quad h\u=l\, M(\x)M(\y)
$$
where $l$ is a complex constant independent of spectral and 
colored parameters.

   So, up to  the solution transformation ${\bf B}$ and ${\bf C}$ one can 
assume 
$$
\au3=\au2=1,\qquad \au8=\au7
$$ 
without losing generality. In the case  
(1.4a), (1.4b), (1.4c) and (1.4d) can be simplified to the following 12 
equations, 

$$
\begin{array}{l}
a_{{5}}(v,\y,\z)+a_{{5}}(u,\x,\y)a_{{1}}(u+v,\x,\z)\\
-a_{{1}}(u,\x,\y)a_{{5}}(u +v,\x,\z)-a_{{7}}(u,\x,\y)a_{{7}}(u+v,\x,\z)
a_{{6}}(v,\y,\z)=0, \\[4mm]
a_{{7}}(u+v,\x,\z)a_{{1}}(v,\y,\z)+a_{{5}}(u,\x,\y)a_{{5}}(u+v,\x,\z)
a_{{7}}(v,\y,\z)\\
-a_{{1}}(u,\x ,\y)a_{{1}}(u+v,\x,\z)a_{{7}}(v,\y,\z)-a_{{7}}(u,\x,\y)
a_{{4}}(v,\y,\z)=0, \\[4mm] 
a_{{6}}(u,\x,\y)+a_{{1}}(u+v,\x,\z)a_{{6}}(v,\y,\z)\\
-a_{{6}}(u+v,\x,\z)a_{{1}}(v,\y,\z)- a_{{5}}(u,\x,\y)a_{{7}}(u+v,\x,\z)
a_{{7}}(v,\y,\z)=0, \\[4mm] 
a_{{6}}(u,\x,\y)a_{{5}}(v,\y,\z)+a_{{1}}(u+v,\x,\z)\\
-a_{{1}}(u,\x,\y)a_{{1}}(v,\y,\z) -a_{{7}}(u,\x,\y)a_{{4} }(u+v,\x,\z)
a_{{7}}(v,\y,\z)=0, \\[4mm] 
a_{{6}}(u,\x,\y)a_{{7}}(u+v,\x,\z)a_{{1}}(v,\y,\z)+ a_{{5}}(u+v,\x,\z)
a_{{7}}(v,\y,\z)\\
-a_{{1}}(u,\x,\y)a_{{7}}(u+v,\x,\z)a_{{5}}(v ,\y,\z)-a_{{7}}(u,\x,\y)
a_{{6}}(u+v,\x,\z)=0, \\[4mm] 
a_{{7}}(u,\x,\y)a_{{1}}(u+v,\x,\z)a_{{1}}(v,\y,\z)+a_{{4}}(u,\x,\y)
a_{{7}}(v,\y,\z)\\
-a_{{1}}(u,\x,\y)a_{{7}}(u+v,\x,\z)-a_{{7}}(u,\x,\y)a_{{6}}(u+v,\x,\z)
a_{{6}}(v,\y,\z)=0
\ea\eqno(1.10)
$$
plus six equations obtained by interchanging  the sub-indices 
1 and 4 as well as 5 and 6 in each of Equations (1.10). We call the six 
equations {\bf counterparts} of (1.10). 

Now we solve the equations obtained by letting $u=0$ and $\y=\x$ in (1.10) with 
respect to the variables $\{a_1(0,\x,\x), a_4(0,\x,\x), a_5(0,\x,\x), 
a_6(0,\x,\x), a_7(0,\x,\x)\}$. It  is easy to prove
that

{\bf Proposition 1.1}\quad For a solution of equations (1.10), 
weight functions satisfy
$$
\ba1
a_1(0,\x,\x)=a_4(0,\x,\x)=1,\\[4mm]
a_5(0,\x,\x)=a_6(0,\x,\x)=a_7(0,\x,\x)=a_8(0,\x,\x)=0.
\ea
\eqno(1.11)
$$
Otherwise, up to the solution transformations A, B, C, D, and E, we have 
two trivial solutions of  the colored Yang-Baxter equation (1.1). The first is 
$$
\begin{array}{l}
\au1=\au4=\au5=\au6= H\u ,\\[3mm]
\au2=\au3=\au7=\au8=1,
\ea\eqno(1.12a)
$$
where  $H\u$ is an arbitrary function of spectral parameter $u$ and 
colored parameters $\x,\y$ and the second 
$$
\begin{array}{l}
\au1=\au4=\au5=-\au6=\dfrac{F(\x)}{F(\y)}exp(u),\\[3mm]
\au1=\au2=1, \qquad             \au7=\au8=i,
\ea\eqno(1.12b)
$$
where $i^2=-1$ and $F(\x)$  is an arbitrary function of colored parameter $\x$.
\vskip5pt

{\bf Definition 1.2}\quad By a gauge solution of the colored Yang-Baxter 
equation (1.1) we mean the solution whose weight functions  
satisfy $\au2=\au3=1$ and \newline 
$\au7=\au8$  and the condition (1.11).
\vskip6pt

(1.11) is called  an initial condition of gauge solutions. The  
condition is simple but very important. It will be quoted again and again 
in finding  gauge solutions of the colored Yang-Baxter equation. For example, 
taking $v=-u,\, \z=\x$ in  (1.10) and their counterparts and 
using the initial condition (1.11), one has the following equations. 
$$
\ba1
a_6(-u,\y,\x)=-a_6(u,\x,\y),\\[4mm]
a_5(-u,\y,\x)=-a_5(u,\x,\y),\\[4mm]
1-\au1 a_1(-u,\y,\x)+\au6 a_5(-u,\y,\x)-\au7 a_7(-u,\y,\x)=0, \\[4mm]
1-\au4 a_4(-u,\y,\x)+\au5 a_6(-u,\y,\x)-\au7 a_7(-u,\y,\x)=0. 
\ea\eqno(1.13)
$$
The unitary condition of a solution of Yang-Baxter equation means 
$$
\check{R}(0,\x,\x)=E,\qquad
\check{R}(u,\x,\y)\,\check{R}(-u,\y,\x)=g(u,\x,\y)\,E,
$$
where $E$ is  the unit matrix and $g(u,\x,\y)$ a scalar function. Hence, it is 
easy to get from (1.13) that 
\vskip5pt

{\bf Proposition 1.3:}\quad For gauge solutions $\check{R}(u,\x,\y)$ 
of the colored Yang-Baxter equation (1.1), the unitary condition is
$$
\check{R}(u,\x,\y)\,\check{R}(-u,\y,\x)=(1-\au5\,\au6)\,E. 
$$ 
\vskip5pt

Differentiating both sides of  all equations in (1.10) and their 
counterparts with respect to the 
variable $v$ and letting $v=0,\,\z=\y$, by virtue of the initial 
condition (1.11)  one  immediately obtains 
$$
\begin{array}{l}
m_5(\y)+\au5\,\MM1-\au1\,\MM5-m_6(\y)\, \au7^{2}=0,\\[4mm]
\MM7+ (m_1(\y)-m_4(\y))\,\au7+m_7(\y)\,(\au5^{2} -\au1^{2})=0,\\[4mm]
\MM6-m_6(\y)\,\au1+m_1(\y)\,\au6 +m_7(\y)\,\au5\,\au7=0,\\[4mm] 
\MM1-m_1(\y)\,\au1+m_5(\y)\,\au6 -m_7(\y)\,\au4\,\au7=0,\\[4mm]
\au6\,\MM7-\au7\, \MM6\\
+m_1(\y)\,\au6\,\au7-m_5(\y)\, \au1\,\au7+m_7(\y)\,\au5=0, \\[4mm] 
\au7\,\MM1-\au1\, \MM7\\
+m_1(\y)\,\au1\,\au7-m_6(\y)\, \au6\,\au7+m_7(\y)\,\au4=0
\end{array}
\eqno(1.14a)
$$
and their counterparts, where and throughout the paper  
we denote 
$$ 
a'_i\u=\dfrac{\partial}{\partial u}a_i\u,
\quad m_i(\x)=a'_i\u \big |_{\{u=0,\y=\x\}}, 
$$ 
for $i=1,4,5,6,7$.  

We call $m_i(\x)$ Hamiltonian coefficients of weight functions with 
respect to spectral parameter or simply coefficients.  Sometimes we write 
$m_i$ instead of $m_i(\x)$ for brevity. 

If we differentiate (1.10)
with respect to $u$ and let $u=0,\,\y=\x$ and then replace the variables 
$v$ and $\z$ by $u$ and $\y$, we have
$$
\begin{array}{l}
m_6(\x)+\au6\,\MM1-\au1\,\MM6-m_5(\x)\, \au7^{2}=0,\\[4mm]
\MM7+ (m_1(\x)-m_4(\x))\,\au7+m_7(\x)\,(\au6^{2} -\au1^{2})=0,\\[4mm]
\MM5-m_5(\x)\,\au1+m_1(\x)\,\au5 +m_7(\x)\,\au6\,\au7=0,\\[4mm] 
\MM1-m_1(\x)\,\au1+m_6(\x)\,\au5 -m_7(\x)\,\au4\,\au7=0,\\[4mm]
\au5\,\MM7-\au7\, \MM5\\
+m_1(\x)\,\au5\,\au7-m_6(\x)\, \au1\,\au7+m_7(\x)\,\au6=0, \\[4mm] 
\au7\,\MM1-\au1\, \MM7\\
+m_1(\x)\,\au1\,\au7-m_5(\x)\, \au5\,\au7+m_7(\x)\,\au4=0
\end{array}
\eqno(1.14b)
$$
and their counterparts.
\vskip5pt

{\bf Remark 1} : In this paper 
, the following trick will often be employed 
to obtain the equations like (1.14a) and (1.14b):
We first do the calculation with respect to $v$ and take $\z=\y$ to 
yield some equations. Then repeat the same operation with $u$ in place of 
$v$ and take $\y=\x$ to yield another equation. 
Then we compare the two results to reduce  to some of the formulas. For example,
formula (1.8) is obtained by this  method.

Comparing (1.14a) with (1.14b), we found that 
the trick, in fact,  is to interchange the subindexes 5 and 6  (or 1 and 4 
) and then replace $m_i(\y)$ by $m_i(\x)$ in original equations (1.14a).
We call this trick {\bf symmetric operation }. 
\vskip5pt

{\bf Remark 2:}\quad If the variable with respect to which we differentiate 
(1.10) is not a spectral parameter $v$ but a colored parameter $\z$, the 
equations 
obtained are the same as (1.14a) if we still let $v=0$ and $\z=\y$. Of course, 
then the dot means the derivative with respect to the colored parameter 
$\y$, that is the second colored variable in $a_i\u$, and 
$m_i(\y)=da_i(v,\y,\z)/d\lambda\big |_{v=0,\z=\y}$. 
Similarly, (1.14b) also
represents the equations obtained by differentiating (1.10) with respect to the 
first colored variable $\x$ in $a_i\u$, but in this case the dot means the 
derivative with respect to $\x$ and
$m_i(\x)=da_i(u,\x,\y)/d\lambda\big |_{u=0,\y=\x}$. 

\vskip10pt 
\subsubsection*{\large \S 2. \bf The coefficients, curves and 
differential equations of weight functions } 

\quad

In this section we will discuss  properties of Hamiltonian coefficients, 
curves and differential
 equations satisfied by weight functions.

It follows from the second equation of (1.14a) and its counterpart that
$$
2a'_7\u+m_7(\y)\,(\au5^2+\au6^2-\au1^2-\au4^2)=0.
\eqno(2.1)
$$
If the symmetric operation is used then one  also has
$$
2a'_7\u+m_7(\x)\,(\au5^2+\au6^2-\au1^2-\au4^2)=0.
\eqno(2.2)
$$
We know $(u_5^2+u_6^2-u_1^2-u_4^2)\neqv 0$ due  to the initial 
condition (1.11). Compared the formulas (2.1) and (2.2), one shows 
$$
m_7(\x)=m_7(\y).
$$
\vskip6pt

If $m_7(\x)=0$ then 
$$
a_7'\u=0,\qquad m_1(\y)=m_4(\y)
$$
due to the second of (1.14a) and its counterpart. This implies $\au7$ is 
a function of only the colored variables $\x$ and $\y$. 
\vskip6pt 

{\bf Proposition 2.1:} \quad For gauge eight-vertex solutions, 
$m_7(\x)$ is a constant independent of colored parameters and $\au7$ 
is independent of spectral parameter if $m_7(\y)=0$.
\vskip5pt

{\bf Remark 3:}\quad In fact, as mentioned in remark 2,  
the latter property of proposition 2.1 also holds for colored parameter, 
namely, $\au7$ will 
be independent of the colored parameter $\x$ (or $\y$) if 
$d/d\x(\au7)\big |_{u=0,\y=\x}=0$ (or $d/d\y(\au7)\big |_{u=0,\x=\y}=0$). 
\vskip5pt

In what follows we denote $m_7(\x)=\alpha$. 
 
Furthermore, if $\alpha=0$, letting $u=0$, $\y=\x$ in the following equation
$$
\MM1=m_6(\y)\,\au6-m_1(\y)\,\au1,
\eqno(2.3)
$$
which is from the  sixth in (1.14a), we obtain
$m_1(\y)=0$ 
owing to the initial condition (1.11). Therefore $m_4(\x)=0$. 
Similarly, we have
$$
\MM6=m_1(\y)\, \au6-m_5(\y)\,\au1 
\eqno(2.4)
$$
from the fifth in (1.14a) and then
$$
m_6(\y)=-m_5(\y).
\eqno(2.5)
$$
Substituting $m_6(\y)=-m_5(\y)$ and $m_1(\y)=m_4(\y)=0$ into
the third and fourth of (1.14a) and their counterparts, we have the following
differential equations. 
$$
\begin{array}{ll}
\MM1=-m_5(\y)\au6 , \qquad & \MM4=m_5(\y)\au5 , \\[3mm]
\MM6=-m_5(\y)\au1 ,        & \MM5=m_5(\y)\au4 .        
\ea\eqno(2.6)
$$
Therefore, if $\alpha=0$ then weight functions satisfy 
$$ \ba1 
\dfrac{d^2}{du^2}a_i(u,\x,\y)=m_5(\y)^2\, a_i\u,\qquad i=1,4,5,6.  
\ea\eqno(2.7) 
$$

\noindent Furthermore, using the symmetric operation  
$$ 
\ba1 \dfrac{d^2}{du^2}a_i(u,\x,\y)=m_5(\x)^2\, a_i\u,\qquad i=1,4,5,6  
\ea\eqno(2.8)
 $$
hold.
Hence $ m_5(\y)$ actually is a constant independent of colored parameters and
not identically zero,
otherwise, the solutions will be independent of spectral parameter. 
Thus we can let $m_5(\y)=\beta$.

The argument above implies  the following proposition.

{\bf Proposition 2.2}\quad For a gauge eight-vertex solution, there exists 
at least 
one between $m_7(\x)$ and $m_5(\x)$ (or $m_6(\x)$) which is not zero 
identically. 
Otherwise, the solution will be independent of spectral parameter.

As for the Hamiltonian coefficients $m_1(\x)$ and $m_4(\x)$ we have  
\vskip5pt

{\bf Proposition 2.3$^{*}$:}\quad For gauge eight-vertex solutions of the colored 
Yang-Baxter equation (1.1)  
$$ m_1(\x)^2-m_4(\x)^2=0.  
$$ 

{\bf Proof:}\quad As first step, we regard the weight functions 
$w_i \,( i=1,4,5,6,7)$
in (1.10) and their counterparts as indeterminates, 
The left side of each of (1.10) and its counterpart  are
polynomial functions of the indeterminates. After eliminating the five indeterminates
$\{w_1,w_4,w_5,w_6,w_7\}$ in numerically increasing order with respect to the 
system of equations (1.10), we can 
obtain seven equations which do not contain the indeterminates 
$w_i\,(i=1,4,5,6,7)$.
 Then  differentiating them with respect to 
the spectral variable $v$ ,  letting $v=0,\z=\y$ and then substituting 
the initial values (1.11) into the resulting ones, we obtain the following 
seven polynomial equations
$$
\ba1
m_{{7}}u_{{1}}^{3}-m_{{7}}u_{{1}}u_{{5}}^{2}-3\,m_{{1}}u_{{1}}u_{{7}}
+m_{{4}}u_{{1}}u_{{7}}-m_{{7}}u_{{4}}u_{{7}}^{2}-m_{{7}}u_{{4}}+m_{{5}}
u_{{6}}u_{{7}}+m_{{6}}u_{{6}}u_{{7}}=0,  \\[4mm]    
-m_{{6}}u_{{1}}u_{{4}}+m_{{7}}u_
{{1}}u_{{6}}u_{{7}}+m_{{7}}u_{{4}}u_{{5}}u_{{7}}+m_{{1}}u_{{4}}u_{{6}}
+m_{{4}}u_{{4}}u_{{6}}-m_{{6}}u_{{5}}u_{{6}}-m_{{5}}u_{{7}}^{2}+m_{{6}}
=0,  \\[4mm]    
m_{{7}}u_{{1}}^{2}-m_{{7}}u_{{4}}^{2}-m_{{7}}u_{{5}}^{2}+m_{{7}}u_{{
6}}^{2}+2\,m_{{4}}u_{{7}}-2\,m_{{1}}u_{{7}}=0,  \\[4mm]    
-m_{{5}}u_{{1}}u_{{4}}+m_{
{1}}u_{{1}}u_{{5}}+m_{{4}}u_{{1}}u_{{5}}+m_{{7}}u_{{1}}u_{{6}}u_{{7}}+
m_{{7}}u_{{4}}u_{{5}}u_{{7}}-m_{{5}}u_{{5}}u_{{6}}-m_{{6}}u_{{7}}^{2}+
m_{{5}}=0,  \\[4mm]    
m_{{7}}u_{{1}}^{2}u_{{6}}-m_{{6}}u_{{1}}u_{{7}}-m_{{5}}u_{{1}}
u_{{7}}-m_{{7}}u_{{5}}^{2}u_{{6}}+m_{{7}}u_{{5}}u_{{7}}^{2}+m_{{7}}u_{
{5}}+m_{{4}}u_{{6}}u_{{7}}+m_{{1}}u_{{6}}u_{{7}}=0,  \\[4mm]    
m_{{7}}u_{{1}}^{3}u_{
{5}}-m_{{6}}u_{{1}}u_{{4}}u_{{7}}-m_{{7}}u_{{1}}u_{{5}}^{3}+2\,m_{{4}}
u_{{1}}u_{{5}}u_{{7}}-2\,m_{{1}}u_{{1}}u_{{5}}u_{{7}}+m_{{7}}u_{{1}}u_
{{6}}  \\
-m_{{7}}u_{{4}}u_{{5}}u_{{7}}^{2}+m_{{5}}u_{{5}}u_{{6}}u_{{7}}+m_
{{6}}u_{{7}}^{3}-m_{{5}}u_{{7}}=0,  \\[4mm]    
-m_{{7}}u_{{1}}^{2}u_{{4}}+m_{{7}}u_{{
1}}u_{{7}}^{2}+m_{{7}}u_{{1}}+m_{{7}}u_{{4}}u_{{5}}^{2}+m_{{4}}u_{{4}}
u_{{7}}+m_{{1}}u_{{4}}u_{{7}}-m_{{5}}u_{{5}}u_{{7}}-m_{{6}}u_{{5}}u_{{
7}}=0.
\ea\eqno(2.9)
$$
As second step, we think of $m_i$ as indeterminates and first
eliminate $m_1,m_4$ and $m_5$ to get two systems of equations, which are 
equivalent to (2.9). The first is
$$
\ba1
-m_{{5}}u_{{1}}u_{{4}}+m_{{1}}u_{{1}}u_{{5}}+m_{{4}}u_{{1}}u_{{5}}+m_
{{7}}u_{{1}}u_{{6}}u_{{7}}+m_{{7}}u_{{4}}u_{{5}}u_{{7}}-m_{{5}}u_{{5}}
u_{{6}}-m_{{6}}u_{{7}}^{2}+m_{{5}}=0,  \\[4mm]    
-m_{{7}}u_{{1}}^{3}u_{{5}}+m_{{7}}u
_{{1}}u_{{4}}^{2}u_{{5}}+2\,m_{{5}}u_{{1}}u_{{4}}u_{{7}}+m_{{7}}u_{{1}
}u_{{5}}^{3}-m_{{7}}u_{{1}}u_{{5}}u_{{6}}^{2}-4\,m_{{4}}u_{{1}}u_{{5}}
u_{{7}}  \\
-2\,m_{{7}}u_{{1}}u_{{6}}u_{{7}}^{2}-2\,m_{{7}}u_{{4}}u_{{5}}u_
{{7}}^{2}+2\,m_{{5}}u_{{5}}u_{{6}}u_{{7}}+2\,m_{{6}}u_{{7}}^{3}-2\,m_{
{5}}u_{{7}}=0,  \\[4mm]    
-m_{{7}}u_{{1}}u_{{4}}^{2}u_{{5}}+m_{{6}}u_{{1}}u_{{4}}u_{
{7}}+m_{{7}}u_{{1}}u_{{5}}u_{{6}}^{2}-m_{{7}}u_{{1}}u_{{6}}+m_{{7}}u_{
{4}}u_{{5}}u_{{7}}^{2}-m_{{5}}u_{{5}}u_{{6}}u_{{7}}  \\
-m_{{6}}u_{{7}}^{3}+m_{{5}}u_{{7}}=0.
\ea\eqno(2.10)
$$
which contains $m_1,m_4$ and $m_5$. The second is
$$
\ba1
\left (u_{{1}}u_{{4}}+u_{{5}}u_{{6}}-u_{{7}}^{2}-1\right )\left (-m_{
{7}}u_{{1}}u_{{4}}^{2}u_{{5}}+m_{{6}}u_{{1}}u_{{4}}u_{{7}}
+m_{{7}}u_{{1}}u_{{5}}u_{{6}}^{2}\right .\\
\left .-m_{{7}}u_{{1}}u_{{6}}+m_{{7}}u_{{4}}u_{{5}}
-m_{{6}}u_{{5}}u_{{6}}u_{{7}}\right )=0,  \\[4mm]   
u_{{1}}\left (u_{{1}}u_{{4}}+u_{{5}}u_{{6}}
-u_{{7}}^{2}-1\right )\left (m_{{7}}u_{{1}}u_{{4}}u_{{5}}^{2}-m_{{6}}u
_{{1}}u_{{5}}u_{{7}}-m_{{7}}u_{{4}}^{2}u_{{5}}u_{{6}}+m_{{6}}u_{{4}}u_
{{6}}u_{{7}}\right .  \\
\left . -m_{{7}}u_{{5}}^{3}u_{{6}}+m_{{7}}u_{{5}}^{2}+m_{{7}}u_{{5
}}u_{{6}}^{3}-m_{{7}}u_{{6}}^{2}\right )=0,  \\[4mm]    
u_{{1}}\left (u_{{1}}u_{{4}}
+u_{{5}}u_{{6}}-u_{{7}}^{2}-1\right )\left (-m_{{7}}u_{{4}}^{3}u_{{5}}
u_{{6}}+m_{{6}}u_{{4}}^{2}u_{{6}}u_{{7}}+m_{{7}}u_{{4}}u_{{5}}^{2}u_{{
7}}^{2}+m_{{7}}u_{{4}}u_{{5}}u_{{6}}^{3}\right .  \\
\left . -m_{{7}}u_{{4}}u_{{6}}^{2}-m_{
{6}}u_{{5}}^{2}u_{{6}}u_{{7}}-m_{{6}}u_{{5}}u_{{7}}^{3}+m_{{6}}u_{{5}}
u_{{7}}\right )=0,  \\[4mm]    
u_{{1}}\left (u_{{1}}u_{{4}}+u_{{5}}u_{{6}}-u_{{7}}^{
2}-1\right )\left (-m_{{7}}u_{{1}}u_{{5}}^{2}u_{{6}}+m_{{7}}u_{{1}}u_{
{5}}-m_{{7}}u_{{4}}^{3}u_{{5}}+m_{{6}}u_{{4}}^{2}u_{{7}}\right.  \\
\left. +m_{{7}}u_{{4}
}u_{{5}}^{3}+m_{{7}}u_{{4}}u_{{5}}u_{{6}}^{2}-m_{{7}}u_{{4}}u_{{6}}-m_
{{6}}u_{{5}}^{2}u_{{7}}\right )=0,
\ea\eqno(2.11)
$$
which do not contain $m_1,m_4$ and $m_5$. So the free fermion 
condition[17]
$$
u_{{1}}u_{{4}}+u_{{5}}u_{{6}}- 1-u_7^2 =0
\eqno(2.12)
$$
or
$$
\ba1
-m_{{7}}u_{{1}}u_{{4}}^{2}u_{{5}}+m_{{6}}u_{{1}}u_{{4}}u_{{7}}+m_{{7}
}u_{{1}}u_{{5}}u_{{6}}^{2}-m_{{7}}u_{{1}}u_{{6}}+m_{{7}}u_{{4}}u_{{5}}
-m_{{6}}u_{{5}}u_{{6}}u_{{7}}=0,  \\[4mm]    
m_{{7}}u_{{1}}u_{{4}}u_{{5}}^{2}-m_{{6}}
u_{{1}}u_{{5}}u_{{7}}-m_{{7}}u_{{4}}^{2}u_{{5}}u_{{6}}+m_{{6}}u_{{4}}u
_{{6}}u_{{7}}-m_{{7}}u_{{5}}^{3}u_{{6}}+m_{{7}}u_{{5}}^{2} \\
+m_{{7}}u_{{5}}u_{{6}}^{3}-m_{{7}}u_{{6}}^{2}=0,  \\[4mm]    
-m_{{7}}u_{{4}}^{3}u_{{5}}u_{{6}}+m_
{{6}}u_{{4}}^{2}u_{{6}}u_{{7}}+m_{{7}}u_{{4}}u_{{5}}^{2}u_{{7}}^{2}+m_
{{7}}u_{{4}}u_{{5}}u_{{6}}^{3}-m_{{7}}u_{{4}}u_{{6}}^{2}-m_{{6}}u_{{5}
}^{2}u_{{6}}u_{{7}}  \\
-m_{{6}}u_{{5}}u_{{7}}^{3}+m_{{6}}u_{{5}}u_{{7}}=0,  \\[4mm]    
-m_{{7}}u_{{1}}u_{{5}}^{2}u_{{6}}+m_{{7}}u_{{1}}u_{{5}} 
-m_{{7}}u_{{4}}^{ 3}u_{{5}}+m_{{6}}u_{{4}}^{2}u_{{7}}
+m_{{7}}u_{{4}}u_{{5}}^{3}+m_{{7}}u _{{4}}u_{{5}}u_{{6}}^{2}
-m_{{7}}u_{{4}}u_{{6}}  \\
-m_{{6}}u_{{5}}^{2}u_{{7}}=0
\ea\eqno(2.13)
$$
will hold.  In the third step, applying the fourth in (2.13) as a main 
equation to kill the 
indeterminate $m_6$
in the three ones remained in (2.13) and then doing factorization of the new 
polynomial equations after killing $m_6$,  one can obtain that  
$$
\ba1
-m_{{7}}u_{{1}}u_{{5}}^{2}u_{{6}}+m_{{7}}u_{{1}}u_{{5}} 
-m_{{7}}u_{{4}}^{ 3}u_{{5}}+m_{{6}}u_{{4}}^{2}u_{{7}}
+m_{{7}}u_{{4}}u_{{5}}^{3}+m_{{7}}u _{{4}}u_{{5}}u_{{6}}^{2}
-m_{{7}}u_{{4}}u_{{6}}\\  
-m_{{6}}u_{{5}}^{2}u_{{7 }}=0,  \\[4mm]  
m_{{7}}\left (u_{{5}}u_{{6}}-1\right )\left (u_{{1}}^{2}u_{{5}}-2\,u_
{{1}}u_{{4}}u_{{6}}+u_{{4}}^{2}u_{{5}}-u_{{5}}^{3}+u_{{5}}u_{{6}}^{2}
\right )u_{{5}}=0,  \\[4mm]   
m_{{7}}\left (u_{{5}}u_{{6}}-1\right )\cdot \\
\left (-u_{{1}}u
_{{4}}^{2}u_{{6}}+u_{{1}}u_{{5}}^{2}u_{{6}}+u_{{1}}u_{{5}}u_{{7}}^{2}-
u_{{1}}u_{{5}}+u_{{4}}^{3}u_{{5}}-u_{{4}}u_{{5}}^{3}-u_{{4}}u_{{6}}u_{
{7}}^{2}+u_{{4}}u_{{6}}\right )u_{{5}}=0,  \\[4mm]    
m_{{7}}\left (u_{{5}}u_{{6}}-1
\right )\left (-u_{{1}}^{2}u_{{4}}+2\,u_{{1}}u_{{5}}u_{{6}}+u_{{4}}^{3
}-u_{{4}}u_{{5}}^{2}-u_{{4}}u_{{6}}^{2}\right )u_{{5}}=0
\ea\eqno(2.14)
$$
is equivalent to (2.13). It follows from the first of (2.14) that  
$m_7\neqv 0$. Otherwise, thanks to proposition 2.2  $\au4^2=\au5^2$. The 
latter  is impossible thanks to the initial condition (1.11).
Hence the following three equations 
$$
\ba1
u_{{1}}^{2}u_{{5}}-2\,u_{{1}}u_{{4}}u_{{6}}+u_{{4}}^{2}u_{{5}}
-u_{{5}}^{3}+u_{{5}}u_{{6}}^{2}=0,  \\[4mm]   
-u_{{1}}u_{{4}}^{2}u_{{6}}+u_{{1}}u_{{5}}^{2}u_{{6}}
+u_{{1}}u_{{5}}u_{{7}}^{2}-u_{{1}}u_{{5}}+u_{{4}}^{3}u_{{5}}
-u_{{4}}u_{{5}}^{3}-u_{{4}}u_{{6}}u_{{7}}^{2}+u_{{4}}u_{{6}}=0,  \\[4mm]    
-u_{{1}}^{2}u_{{4}}+2\,u_{{1}}u_{{5}}u_{{6}}+u_{{4}}^{3}
-u_{{4}}u_{{5}}^{2}-u_{{4}}u_{{6}}^{2}=0
\ea\eqno(2.15)
$$
and the first of (2.14) is equivalent to (2.14), where we use the initial 
condition (1.11) again to yield $1-u_5u_6\neqv 0$.

Finally, if  we differentiate (2.12) and the third equation of (2.15) 
with respect to $u$ and let $u=0$, $\x=\y$ and then
apply the initial conditions (1.11) again , we come to the conclusion of 
proposition 2.3.
\vskip5pt

{\bf Remark 4:}\quad When we do the operation of eliminating indeterminates
with respect to a system of polynomial equations, according to the 
theorem of zero 
structure of algebraic varieties [22],  the coefficient of  the term with 
 the highest 
degree of the indeterminate in main polynomial equation (to be  
eliminated in other 
polynomials) should not be  identified with zero. In the event it is 
identified with 
zero, we should add the coefficient into the system of equations to produce
a new system of equations. Otherwise, it is possible to lose some solutions.
For example, when we use the fourth equation in (2.13) to eliminate $m_6$ 
in the three remaining ones in (2.13), because the coefficient of $m_6$ in 
the fourth one of (2.13) is $u_7(u_4^2-u_5^2)$ which does not identify with 
zero 
due to the initial condition (1.11), (2.14) is equivalent to the system
of equations (2.13).
\vskip6pt

From the argument of proving proposition 2.3 above we see the system of 
equations
 (2.9) is equivalent to two groups of equations. The first is (2.10),
(2.15) plus the first of (2.14). The second is (2.10) and (2.12), 
the free fermion condition.

Now we consider the two cases respectively. For the first case we differentiate 
the second equation in (2.15) and take $u=0,\, \x=\y$. Then we 
substitute the initial condition (1.11) into the result to get 
$$ 
m_5(\y)=m_6(\y).  
\eqno(2.16) 
$$ 
By doing factorization of the equation obtained 
by eliminating  $u_4$ in the third equation of (2.15)
and by using the second equation of (2.15) , we can get 
$$ 
2\,u_{{6}}\left (u_{{6}}-u_{{5}}\right )\left (u_{{6}}+u_{{5}}\right ) 
\left (u_{{1}}-u_{{5}}\right )\left (u_{{1}}+u_{{5}}\right )
\left (u_{ {1}}-u_{{6}}\right )\left (u_{{1}}+u_{{6}}\right )=0. 
\eqno(2.17) 
$$ 
Together, (2.17) and the initial condition (1.11) imply 
$$
a_5(u,\x,\y)=a_6(u,\x,\y) 
\eqno(2.18a) 
$$ 
or
$$
a_5(u,\x,\y)=-a_6(u,\x,\y) 
\eqno(2.18b) 
$$
hold. Substituting (2.18a) into the third equation in (2.15), we then see 
$$ 
a_1(u,\x,\y)=a_4(u,\x,\y).  
\eqno(2.19) 
$$ 
Substituting (2.18a) and (2.19) into the first equation in (2.14), we have 
$$ 
(u_5- u_1 (u_5+u_1 )( \alpha\,u_5 u_1-m_6\, u_7 )=0.  
\eqno(2.20) 
$$ 
So 
$$
\alpha\,\au1\,\au5 -m_6(\y)\,\au7 =0,
\eqno(2.21)
$$
here the initial condition (1.11) is used again. (2.21) implies 
$m_6(\y)\neqv 0$, or $\alpha$ will also identify with zero. But this will 
contradict proposition 2.2. 
 If we do  the symmetric operation with respect to (2.21) and use (2.19), 
then 
$$ 
\alpha\,\au1\,\au5 -m_6(\x)\,\au7 =0,
\eqno(2.22) 
$$ (2.21) and (2.22) will imply $m_6(\y)$  is also a constant independent 
of colored parameter. We let $m_5(\y)=\beta$.  

If $\au5=-\au6$ we should have $m_6(\y)=0$ and (2.21) still holds. Then 
$\alpha=0$.  It is clearly impossible thanks to proposition 2.2.

Combining (2.18a), (2.19), (2.21) with  the third equation 
in (1.14a) we obtain  
$$
(a_{{5}}'\u)^2 =\beta^2- (\beta^2-m_1(\y)^2+\alpha^2)\,\au5^2 
+\alpha^2\,\au5^4.
\eqno(2.23) 
$$
Using the symmetric operation we can show $m_1(\y)$ is also a constant 
independent of colored parameter. Let $m_1(\y)=\gamma$.
Similarly,
$$ 
(a_{{1}}'\u)^2 
=\beta^2- (\beta^2-\gamma^2+\alpha^2)\,\au1^2 +\alpha^2\,\au1^4 
\eqno(2.24) 
$$
holds if we combine (2.18a), (2.19), (2.21) with the fourth in (1.14a). 
Substituting (2.16),\newline 
(2.18a), (2.19) and (2.21) into the second of (2.10) 
we fund that the algebraic curve satisfied by the weight functions 
$\au1$ and $\au5$ is 
$$ 
\alpha^2\,u_1^2u_5^2-\beta^2\,u_5^2-\beta^2u_1^2
+2\beta\,\gamma\, u_1u_5+\beta^2=0. 
\eqno(2.25)
$$

From the second group of equations, i.e. (2.10) and (2.12),  
it is easy to obtain
$$
\ba1
1+\au7^2-\au1\au4-\au5\au6=0,\\[4mm]
\alpha\, (\au1\au6+\au4\au5)=(m_5(\y)+m_6(\y))\,\au7,\\[4mm]
\alpha\, (\au1^2+\au6^2-\au4^2-\au5^2)=4m_1(\y)\,\au7.
\ea
\eqno(2.26)
$$
Then $m_1(\y)+m_4(\y)=0$.
If we do the symmetric operation with respect to the third equation of (2.26) 
then we also obtain 
$$ 
\alpha\, (\au1^2+\au5^2-\au4^2-\au6^2)=4m_1(\x)\,\au7.  
\eqno(2.27) 
$$ 
Letting $\y=\x$ in the  third equation in (2.26) and (2.27) we reduce to 
$$ 
a_6(u,\x,\x)^2=a_5(u,\x,\x)^2.  
\eqno(2.28)
$$
Therefore, (2.16) and (2.28) give the following proposition. 
\vskip5pt
 
{\bf Proposition 2.4:} \quad For gauge eight-vertex solutions of the colored
 Yang-Baxter equation (1.1) the Hamiltonian coefficients $m_5(\y)$ 
and $m_6(\y)$ satisfy
$$
m_5(\y)^2-m_6(\y)^2=0.
$$

From (2.26) and the second equation in (1.14a), we can calculate the 
weight function $\au7$ obeys 
$$ 
\ba1 
\qquad a_7'(u,\x,\y)^2\\ 
=\alpha^2-((m_{{5}}(\y)+m_6(\y))^2-4\,m_{{1}}^{2}(\y)-2\,\alpha^{2}) 
\au7^{2}+\alpha^2\au7^4.  
\ea\eqno(2.29) 
$$ 
Furthermore, as we said in remark 1, we can show 
$(m_5(\y)^2+m_6\y))^2-2m_1(\y)^2$ 
is a constant independent  of colored 
parameter using the symmetric operation. Let
$\delta^2=(m_5(\y)^2+m_6(\y))^2-2m_1(\y)^2$. 

We conclude this section by the following theorem.
\vskip5pt

{\bf THEOREM 2.5$^{*}$:}\quad For a gauge eight-vertex type solution, its 
weight functions must satisfy one of  two systems of equations. The first 
is composed of (2.18a), (2.19), (2.21), (2.23), (2.24) and (2.25).  
The second is composed of (2.26) and (2.29).
\vskip10pt

\subsubsection*{\large \S 3. \bf  Gauge eight-vertex type solutions } 

\quad

In this section $sn(\zeta)$ and $cd(\zeta)=cn(\zeta)/dn(\zeta)$ are 
Jacobian elliptic functions. 

Now  we describe how to write down  all gauge solutions of eight-vertex 
type of the colored Yang-Baxter (1.1) and classify them into two types called 
Baxter type and Free-Fermion type.  
\vskip5pt 

{\bf 3.1  Baxter type solutions.} 

\quad

We  consider the first case in theorem 2.6.
Since the curve (2.25) only includes two weight functions $\au1$ and $\au5$ 
, we can parameterize 
$\au1$ and $\au5$ as one-parameter functions. If $\beta\pm \alpha\pm \gamma\not= 0$, 
(2.23) and (2.24) are elliptic differential equations. Therefore,  
the  solutions should be 
$$
\begin{array}{l}
\au1=\au4=\dfrac{sn(\lambda\,u+F(\x)-F(\y)+\mu)}{sn{(\mu)}},\\[4mm]
\au2=\au3=1,\\[4mm]
\au5=\au6=\pm\dfrac{sn(\lambda\,u+F(\x)-F(\y))}{sn{(\mu)}},\\[4mm]
\au7=\pm k\,sn(\lambda\,u+F(\x)-F(\y))\,sn(\lambda\,u+F(\x)-F(\y)+\mu)
\end{array}
\eqno(3.1)
$$
where $k$, as the modules of  Jacobi elliptic function, is an arbitrary constant.

If
 $\beta\pm \alpha\pm \gamma=0$, the elliptic solutions (3.1) will degenerate 
into  trigonometric solutions
$$
\begin{array}{l}
\au1=\au4=\dfrac{tan(\lambda\,u+F(\x)-F(\y)+\mu)}{tan(\mu)},\\[4mm]
\au2=\au3=1,\\[4mm]
\au5=\au6=\pm\dfrac{tan(\lambda\,u+F(\x)-F(\y))}{tan(\mu)},\\[4mm]
\au7=\pm tan(\lambda\,u+F(\x)-F(\y))\,tan(\lambda\,u+F(\x)-F(\y)+\mu). 
\ea\eqno(3.2)
$$
In (3.1) and (3.2) $\lambda\not= 0$, $\mu \not= 0$ are two arbitrary constants  
and $F(\x)$ an arbitrary function. 
\vskip7pt

{\bf 3.2   Free-Fermion type solutions.} 

\quad

Now we consider the second case in theorem (2.6). According to 
proposition 2.4 it is divided into two sub-cases, 
$m_5(\x)=m_6(\x)$ and $m_5(\x)=-m_6(\x)$.
\vskip5pt

{\bf (3.2a)   For the sub-case of $m_5(\x)=m_6(\x)$} 

Let $m_5(\x)\neqv 0$ (we will put the case of $m_5(\y)=0$ into the second 
sub-case ). It is clear that $\alpha\not= 0$ due to the second of (2.26). 
For brevity  we let $\alpha=1$ up to the solution transformation D. 
When 
$$
m_5(\x)^2-m_1(\x)^2\neqv 0 
$$
and 
$$
m_5(\x)^2-m_1(\x)^2\neqv 1 
$$
the equation (2.29) then is  an elliptic  differential equation 
and, from the remarks 2, 3 and the initial condition (1.11),  should have 
solutions
$$
\au7=k\,sn(\lambda\,u+F(\x)-F(\y))\,cd(\lambda\,u+F(\x)-F(\y)),
\eqno(3.3)
$$
where $k$, as 
the module of elliptic function, and  $\lambda$ are two arbitrary constants 
and $F(\x)$ an arbitrary function with the constriction $k\,\lambda=1$. 

Substituting (3.3) into the second of (1.14a) and its counterpart 
as well as the first and second of (2.26) and using elliptic function 
identities  we have
\renewcommand{\theequation}{3.4\alph{equation}} 
\setcounter{equation}{0} 
\begin{eqnarray} 
{\it cd}^{2}-{\it sn}^{2}+2\,m_{{1}}(\y)k\,{\it cd}\,{\it sn}+u_{{5}}^{2}
-u_{{1}}^{2}=0,\\[4mm]
{\it cd}^{2}-{\it sn}^{2}-2\,m_{{1}}(\y)k\,{\it cd}\,{\it sn}+u_{{6}}^{2}
-u_{{4}}^{2}=0,\\[4mm]
u_{{1}}u_{{4}}+u_{{5}}u_{{6}}-{\it cd}^{2}-{\it sn}^{2}=0,\\[4mm]
u_{{1}}u_{{6}}+u_{{4}}u_{{5}}-2\,m_{{5}}(\y)k\,{\it cd}\,{\it sn}=0,
\end{eqnarray} 
where $m_5(\y)$ and $m_1(\y)$ are arbitrary functions satisfying 
$$
m_5(\y)^2-m_1(\y)^2=\dfrac{1}{k^2}.
$$
In the formulas (3.4)s and in what follows we simply write $sn,\, cd$ instead of 
elliptic functions $\sn$ and $\cd$ for brevity. Using the symmetric 
operation we also have 
$$
{\it cd}^{2}-{\it sn}^{2}+2\,m_{{4}}(\x)k\,{\it cd}\,{\it sn}+u_{{5}}^{2}
 -u_{{4}}^{2}=0. 
$$
Since $m_1(\x)+m_4(\x)=0$ one obtains  
$$
{\it cd}^{2}-{\it sn}^{2}-2\,m_{{1}}(\x)k\,{\it cd}\,{\it sn}+u_{{5}}^{2} 
-u_{{4}}^{2}=0.  
\eqno(3.4e)
$$
From (3.4c), (3.4d) and (3.4a) one also obtains
$$
-({\it cd}^{2}+{\it sn}^{2})u_{{1}}+ ({\it cd}^{2}-{ \it sn}^{2}
+2\,m_{{1}}(\y)k\,{\it sn}\,{\it cd})u_{{4}}+2\,m_{{5}}(\y)k\,{\it sn}
\,{\it cd}\,u_5=0.
\eqno(3.5a)
$$
If we do the symmetric operation with respect to the counterpart of (3.5a) then 
$$
-({\it cd}^{2}+{\it sn}^{2})u_{{4}}+({\it cd}^{2}-{ \it sn}^{2}
-2\,m_{{1}}(\x)k\,{\it sn}\,{\it cd})u_{{1}}+2\, m_{{5}}(\x)k\,
{\it sn}\,{\it cd}\, u_5=0.
\eqno(3.5b)
$$
Similarly, one has
\renewcommand{\theequation}{3.6\alph{equation}} \setcounter{equation}{0}
 \begin{eqnarray} 
({\it sn}^{2}+2\,m_{{1}}(\y)k\,{\it sn}\,{\it cd}-{\it cd}^{2}
)u_{{5}}-({\it cd}^{2}+{\it sn}^{2})u_{{6}}+2\,m_{{5}}(\y)k\,
{\it sn}\,{\it cd}\, u_4=0,\\[4mm]
({\it sn}^{2}+2\,m_{{1}}(\x)k\,{\it cd}\,{\it sn}-{\it cd}^{2}
)u_{{6}}+({\it cd}^{2}+{\it sn}^{2})u_{{5}}+2\,m_{{5}}(\x)k\,{\it cd}\,
{\it sn}\, u_4=0.
\end{eqnarray} 
Solving the equations (3.5a), (3.5b), (3.6a) and (3.6b) with respect to 
$\{u_1,u_4,u_5,u_6\}$ we have
$$
\dfrac{\au4}{\au1}=\dfrac{H_4}{H_1},\qquad\qquad 
\dfrac{\au6}{\au5}=\dfrac{H_6}{H_5},
$$
where
$$
\ba1
H_1=(m_5(\x)+m_5(\y))\,cd+(m_1(\x)m_5(\y)+m_5(\x)m_1(\y))\,k\,sn,\\[4mm] 
H_4=(m_5(\x)+m_5(\y))\,cd-(m_1(\x)m_5(\y)+m_5(\x)m_1(\y))\,k\,sn,\\[4mm] 
H_5=(m_5(\x)+m_5(\y))\,sn+(m_1(\x)m_5(\y)-m_5(\x)m_1(\y))\,k\,cd,\\[4mm]  
H_6=(m_5(\x)+m_5(\y))\,sn-(m_1(\x)m_5(\y)-m_5(\x)m_1(\y))\,k\,cd. 
\ea\eqno(3.6)
$$
Let $u_1=H_1X,\, u_4=H_4X$ and $u_5=H_5Y,\, u_6=H_6Y$. From (3.4a) and 
(3.4e) one obtains
$$
\begin{array}{ll}
(H_1^2-H_4^2)X^2
&=4\,\left (m_{{5}}(\x)+m_{{5}}(\y)\right )\left (m_{{5}}(\y)m_{{1}}(\x)
+m_{{5}}(\x)m_{{1}}(\y)\right )\,k\,{\it sn}\,{\it cd}\,X^2\\[3mm]
&=2(m_1(\x)+m_1(\y)k\,sn\,cd,  \\[4mm]
(H_5^2-H_6^2)Y^2
&= 4\,\left (m_{{5}}(\x)+m_{{5}}(\y)\right )\left (m_{5}(\y)m_{{1}}(\x)
-m_{{5}}(\x)m_{{1}}(\y)\right )k\, {\it sn}\,{\it cd}\,Y^2\\[3mm]
&=2(m_1(\x)-m_1(\y))k\,sn\,cd.
\ea\eqno(3.7)
$$
It is easy to show using $m_5(\x)^2-m_1(\x)^2=1/k^2$ that 
$$
\begin{array}{ll}
X^2&=\dfrac{m_1(\x)+m_1(\y)}{2(m_5(\x)+m_5(\y))(m_1(\x)m_5(\y)+m_5(x)m_1(\y))}
    \\[6mm]
   &=\dfrac{1+k^2(m_5(\x)m_5(\y)- m_1(\x)m_1(\y))}{2(m_5(\x)+m_5(\y))^2}
   =\dfrac{-1+k^2(m_5(\x)m_5(\y)+ m_1(\x)m_1(\y))}
      {2(m_1(\x)m_5(\y)+m_5(\x)m_1(\y))^2\,k^2},\\[7mm]
Y^2&=\dfrac{m_1(\x)-m_1(\y)}{2(m_5(\x)+m_5(\y))(m_1(\x)m_5(\y)-m_5(x)m_1(\y))}
   \\[6mm]
   &=\dfrac{1+k^2(m_5(\x)m_5(\y)+ m_1(\x)m_1(\y))}{2(m_5(\x)+m_5(\y))^2}
   =\dfrac{-1+k^2(m_5(\x)m_5(\y)- m_1(\x)m_1(\y))}
      {2(m_1(\x)m_5(\y)-m_5(\x)m_1(\y))^2\,k^2}.
\ea
$$
Hence  gauge solutions of the colored Yang-Baxter 
equation (1.1) should obey the following forms 
$$
\ba1
\au1= A(\x,\y)\,\cd+ B(\x,\y)\,\sn ,\\[4mm] 
\au2=\au3=1,\\[4mm]
\au4= A(\x,\y)\,\cd- B(\x,\y)\,\sn ,\\[4mm] 
\au5=C(\x,\y)\,\sn+ D(\x,\y)\,\cd  ,\\[4mm] 
\au6=C(\x,\y)\,\sn- D(\x,\y)\,\cd  ,\\[4mm] 
\au7=\pm k\,\sn\,\cd ,
\ea\eqno(3.8)
$$
where  $k$, as  the module of the elliptic functions,  is an arbitrary 
constant and
$$ 
\ba1 
A(\x,\y)= \sqrt{(1+G(\x)G(\y)- H(\x)H(\y))/2},\\[4mm]
B(\x,\y)= \sqrt{(-1+G(\x)G(\y)+ H(\x)H(\y))/2},\\[4mm]
C(\x,\y)= \delta\sqrt{(1+G(\x)G(\y)+ H(\x)H(\y))/2},\\[4mm]
D(\x,\y)= \delta \sqrt{(-1+G(\x)G(\y)- H(\x)H(\y))/2}\,M,
\ea\eqno(3.9) 
$$
where $\delta^2=1$ and
$$
M=\dfrac{H(\x)G(\y)- G(\x)H(\y)}{\sqrt{(H(\x)G(\y)- G(\x)H(\y))^2}}
$$
and $G(\x)$,\, $H(\x)$ satisfy $G(\x)^2-H(\x)^2=1$. If we consider the 
solution transformation D then the restrictive condition 
$k\,\lambda=1$ can be cancelled, namely, $\lambda\not= 0$ is also an arbitrary 
constant.
     
    When the module $k=1$ the Jacobian elliptic functions $cd$ and $sn$ 
should degenerate into $1$ and $tanh$. Hence,  we have 
$$
\ba1
\au1= A(\x,\y)+ B(\x,\y)\,\tanh ,\\[4mm] 
\au2=\au3=1,\\[4mm]
\au4= A(\x,\y)- B(\x,\y)\,\tanh ,\\[4mm] 
\au5=C(\x,\y)\,\tanh+ D(\x,\y),\\[4mm] 
\au6=C(\x,\y)\,\tanh- D(\x,\y),\\[4mm] 
\au7=\pm tanh\ll, 
\ea\eqno(3.10)
$$
where $A(\x,\y),\,B(\x,\y)$, $C(\x,\y)$ and $D(\x,\y)$ are defined by (3.9) 
and  $\lambda\not= 0$ is an arbitrary constant.
\vskip5pt

If $m_5(\y)^2-m_1(\y)^2=0$ then the differential equation (2.29) can be 
rewritten  as
$$
a_7'(u,\x,\y)^2 =\alpha^2(1+2\, \au7^{2}+\au7^4).
$$
So, following the calculation of (3.8) we can show that 
the gauge solutions are 
$$
\ba1
\au1=X\left (\dfrac{ G(\x)+G(\y)}{cos\ll}+2 G(\x)G(\y) 
sin\ll\right ),\\[4mm]
\au2=\au3=1,\\[4mm]
\au4=X\left (\dfrac{ G(\x)+G(\y) }{cos\ll}-2 G(\x)G(\y)
 sin\ll\right ),\\[4mm]
\au5=Y\left (\dfrac{ G(\x)-G(\y) }{cos\ll}+2 G(\x)G(\y) 
sin\ll\right),\\[4mm]
\au6=Y\left (-\dfrac{ G(\x)-G(\y) }{cos\ll}+2 G(\x)G(\y) 
sin\ll \right),\\[4mm]
\au7=\pm tan\ll,
\ea\eqno(3.11)
$$
where
$$
X=\dfrac{1}{2\,\sqrt{G(\x)G(\y)}},\qquad Y=\pm X 
\eqno(3.12)
$$
and $G(\x)$ is an arbitrary function.
\vskip5pt

{\bf (3.2b)   For the sub-case of $m_5(\x)=-m_6(\x)$} 

In the case of $m_5(\y)=-m_6(\y)$ weight functions of gauge solutions 
are 
$$
\ba1
\au1=\au4=\dfrac{cosh(\lambda\,u+F(\x)-F(\y))}{cos(\mu\,u+G(\x)-G(\y))},\\[4mm]
\au2=\au3=1,\\[4mm]
\au5=-\au6=\pm \dfrac{sinh(\lambda\,u+F(\x)-F(\y))}{cos(\mu\,u+G(\x)-G(\y))},
   \\[4mm]
\au7=\pm tan(\mu\,u+G(\x)-G(\y)),
\ea
\eqno(3.13)
$$
where $\lambda$ and $\mu$ are two arbitrary constants, but not zero 
simultaneously,  and $F(\x)$, 
$G(\x)$ are  two arbitrary functions.

To prove it we first consider the sub-case of $\alpha\not= 0$. Then it follows  that
$$
\au1\au6+\au4\au5=0
$$
by (2.26) and
$$
\au1\au5+\au4\au6=0
$$
by the symmetric operation.  Since $\au1\not= -\au4$ owing to the initial 
condition (1.11). So one can get
$$
\ba1
\au1=\au4,\\[4mm]
\au5=-\au6,\\[4mm]
\au1^2-\au5^2=1+\au7^2
\ea\eqno(3.14)
$$ 
and then $m_1(\y)=0$ because of the first of (3.14) and the condition 
$m_1(\x)+m_4(\x)=0$. Then (2.29) has solution 
$$
\au7=tan(\mu\,u+G(\x)-G(\y))
$$
and hence (3.13) is true.

Second we consider the sub-case of $\alpha=0$. Then it follows from the 
argument of proposition 2.2 in section two that 
the weight functions are 
$$ 
\begin{array}{l} 
\au1=A_1(\x,\y)cosh\,u - A_6(\x,\y)sinh\,u, \\[4mm] 
\au4=A_4(\x,\y)cosh\,u + A_5(\x,\y)sinh\,u, \\[4mm] 
\au5=A_5(\x,\y)cosh\,u + A_4(\x,\y)sinh\,u,\\[4mm] 
\au6=A_6(\x,\y)cosh\,u - A_1(\x,\y)sinh\,u,\\[4mm] 
\au7=A_7(\x,\y),  
\end{array}\eqno (3.15) 
$$ 
where $A_i(\x,\y)=a_i(0,\x,\y)\, (i=1,4,5,6,7)$ are some functions with 
respect to colored parameters $\x,\y$ to be determinate.  
It is clear that, 
$A_i(\x,\y)\, (i=1,4,5,6,7)$  satisfy the pure the colored Yang-Baxter (1.2).

If we substituting the initial condition (1.11) into the ones obtained by 
putting $v=-u$ and $\z=\y$ in the system of (1.10) and combining resulting 
equation with the free fermion condition (2.12), one 
can show
$$
\ba1
a_{{5}}(u,\x,\y)=-a_{{5}}(-u,\y,\x),\qquad 
     a_{{6}}(u,\x,\y)=-a_{{6}}(-u,\y,\x) ,\\[4mm] 
\au7=-a_7(-u,\y,\x),\qquad \au4=a_1(-u,\y,\x).  
\ea\eqno(3.16)
$$ 
It is easy by (3.16) to show that
$$
A_4(\x,\y)=A_1(\x,\y),\qquad A_6(\x,\y)=-A_5(\x,\y).
\eqno(3.17)
$$
   
As mentioned in remark 2, all formulas obtained in this section and first two 
sections should hold for the pure the colored Yang-Baxter equation 
(1.2) except those  obtained by using symmetric operation. So
  we still should have 
$$
\ba1
1+A_7(\x,\y)^2-A_1(\x,\y)A_4(\x,\y)-A_5(\x,\y)A_6(\x,\y)=0,\\[4mm]
l_7(\y)(A_1(\x,\y)A_6(\x,\y)+A_4(\x,\y)A_5(\x,\y))=(l_5(\y)+l_6(\y))A_7(\x,\y),
      \\[4mm]
l_7(\y)(A_1(\x,\y)^2+A_6(\x,\y)^2-A_4(\x,\y)^2-A_5(\x,\y)^2
=4\,l_1(\y)A_7(\x,\y),
\ea\eqno(3.18)
$$
where, as mentioned in remark 2, $l_i(\y),\, (i=1,4,5,6,7)$ mean 
$\dfrac{\partial}{\partial \y}A_i(\x,\y)\big|_{\x=\y}$. Substituting (3.17) 
into (3.18) we see 
$$
l_5(\y)+l_6(\y)=0,\qquad l_1(\y)=0.
$$
Therefore, as we did for getting (2.29) we also have 
$$
(\dfrac{\partial}{\partial \y}A_7(\x,\y))^2
 =l_7(\y)(1+2\,\au7^{2}+\au7^4)
\eqno(3.19)
$$
and then the solution for $A_i(\x,\y)$ is
$$
\ba1
A_1(\x,\y)=A_4(\x,\y)=\dfrac{cosh(F(\x)-F(\y))}{cos(G(\x)-G(\y))},\\[4mm] 
A_5(\x,\y)=-A_6(\x,\y)=\dfrac{sinh(F(\x)-F(\y))}{cos(G(\x)-G(\y))},\\[4mm] 
A_7(\x,\y)=tan(G(\x)-G(\y)), 
\ea\eqno(3.20)
$$
where $F$ and $G$ are two arbitrary functions of single variable. 
Substituting (3.20) into (3.15) 
one can say (3.13) is also true for the case of $\alpha=0$, only $\mu=0$.

Description above tell us that for the case of Hamiltonian coefficients 
$m_5(\y)=-m_6(\y)$  weight functions of a gauge eight-vertex 
type solution must be (3.13).
 
Finally, straightforward calculation and computer symbolic computation 
can verify the following theorem.
 \vskip6pt

{\bf THEOREM 3.1 :}\quad  
Gauge eight-vertex solutions of the colored Yang-Baxter equation (1.1) are 
composed of  (3.1), (3.8) and (3.13) and their degenerate 
forms (3.2), (3.10) and (3.11).
\vskip4pt

 The solutions (3.1) and (3.2) are called 
Baxter type solutions. They are just the solutions for "zero field" 
eight-vertex model by Baxter [3]. 
 The solutions (3.8), (3.13) and their degenerate forms (3.10), (3.11) satisfy 
free fermion condition and are called  
Free-Fermion type solutions.  
\vskip6pt

If we take $\lambda=1, \, G(\x)=cosh(2\x),\, 
H(\x)=sinh(2\x)$ and 
$F(\x)=0$ in solution (3.8) then (3.8) will reduce to the following solution 
$$
\ba1
\au1=cosh(\x-\y)cd(u)+sinh(\x+\y)sn(u),\\[4mm] 
\au4=cosh(\x-\y)cd(u)-sinh(\x+\y)sn(u),\\[4mm] 
\au5=cosh(\x+\y)sn(u)-sinh(\x-\y)cd(u),\\[4mm]  
\au6=cosh(\x+\y)sn(u)+sinh(\x-\y)cd(u),\\[4mm] 
\au7=k\,sn(u)\,cd(u)
\ea
$$
which is given in Ref. 20.
\vskip10pt

\subsubsection*{\large \bf \S 4. General solutions}

\quad

In this paper we have shown and classified all gauge eight-vertex solutions of 
the colored Yang-Baxter equation (1.1). These gauge solutions and trivial 
solutions (1.12a), (1.12b) together with five solution transformations 
discussed in the first section will give all eight-vertex type solutions. 

  If  we take in (3.8) and (3.9)
$$
G(\x)=\dfrac{1}{sn(\x)},\qquad 
H(\x)=\dfrac{ cn(\x)}{sn(\x)},
\qquad F(\x)=0,\qquad \lambda=\dfrac 12
$$
 and the solution transformation B with 
$$
g(u,\x,\y)=\sqrt{e(\x)e(\y)sn(\x)sn(\y)}\,\dfrac{(1-e(u))}{sn(u/2)},
$$
where the elliptic exponential
$$
e(\zeta)=cn(\zeta)+i\, sn(\zeta)
$$
then 
using addition theorems for elliptic functions $sn(\zeta), cn(\zeta) $ and 
$dn(\zeta)$ we can obtain the following solution given in Ref. 12
$$
\ba1
\au1=1-e(u)\,e(\x)\,e(\y),\\[4mm]
\au2=\au3=\sqrt{e(\x)e(\y)sn(\x)sn(\y)}\dfrac{(1-e(u))}{sn(u/2)},\\[4mm]
\au4=e(u)-e(\x)\,e(\y),\\[4mm]
\au5=e(\x)-e(u)\,e(\y),\\[4mm]
\au6=e(\y)-e(u)\,e(\x),\\[4mm]
\au7=\au8=-ik\, \sqrt{e(\x)e(\y)sn(\x)sn(\y)}\,(1-e(u))\,sn(u/2),
\ea
$$
the detail of calculation of which is omitted.

Similarly, all non-trivial general solutions can be also classified into 
two types. 
The first are Baxter type solutions if  they can be obtained via gauge Baxter 
solutions and some solution transformations. The second are Free-Fermion 
solutions if they can be obtained via gauge Free-Fermion solutions and some 
solution transformations.

According to the standard method by Baxter, for a given R-matrix the 
spin-chain Hamiltonian is generally of the following form:
$$
H=\sum_{j=1}^N(J_x\,
\sigma_j^x\sigma_{j+1}^x+ J_y\,\sigma_j^y\sigma_{j+1}^y+ 
       J_z\,\sigma_j^z\sigma_{j+1}^z+ \dfrac 12 h\,(\sigma_j^z+\sigma_{j+1}^z)),
$$
where $\sigma^x,\, \sigma^y,\,\sigma^z$ are Pauli matrices and the 
coupling constants are
$$
\ba1
J_x=\dfrac 14 (m_5+m_6+m_7+m_8),\quad  
J_y=\dfrac 14 (m_5+m_6-m_7-m_8), \\[4mm] 
J_z=\dfrac 14 (m_1-m_3+m_4-m_2),\quad  
h=\dfrac 14 (m_1-m_3-m_4+m_2).  
\ea
$$

In this paper we have  proved that the hamiltonian 
coefficients  of a gauge solution must obey 
$$
m_1^2=m_4^2,\qquad m_5^2=m_6^2.
$$
It follows from the solution transformations $ B$ and $D$ that 
$$
(m_1-m_3)^2=(m_4-m_2)^2,\qquad m_5^2=m_6^2
$$
for general solutions. This clearly describes the relation between 
classifications of eight-vertex 
type solutions and spin-chain Hamiltonians. For example, if 
$J_x+J_y=h,\, J_z=0$, i.e. a special free-fermion model in a magnetic 
field, then one has $m_5=m_1-m_3=-m_4+m_2$. The corresponding solution 
of the colored Yang-Baxter equation should be (3.11). Then the  transfer 
matrix is of  
trigonometric function type.
\vskip15pt

\noindent{\bf Acknowledgement:} \quad 
This work was carried out while 
the author was visiting Department of Pure Mathematics, University of Adelaide. 
The author would like to thank Alan Carey for his invitation to the 
department and  the department for kind hospitality.

This work was supported by Climbing Up Project, NSCC and Natural 
Scientific Foundation of Chinese Academy of Sciences.
\vskip 20pt

\subsubsection*{\bf References}

\begin{enumerate}

\item[{[1]}]  Yang, C.N., 
Phys. Rev. Lett., {\bf 19} (1967), 1312-1314.

\item[{[2]}]  Yang, C.N.,
Phys. Rev. Lett., {\bf 168} (1968), 1920-1923.

\item[{[3]}]  Baxter, R.J., 
Ann. Phys. {\bf 70} (1972), 193-228.

\item[{[4]}] Zamolodchikov, A. B. and Zamolodchikov A. B.,
Annals of Physics, {\bf 120} (1979), 253-291.

\item[{[5]}] Baxter, R.J., Exactly solved models in statistical mechanics, 
Academic Press, London, 1982.

\item[{[6]}] Jimbo, M., Yang-Baxter equation in integrable systems,
World Scientific, 1989. 

\item[{[7]}]  Drinfel'd, V.G.,
"Quantum Groups", Proceeding of th International Congress of Mathematicians, 
Berkeley, (1987), 798-820.

\item[{[8]}]  Alvarez-Gaum\'e, L., C\'omez, C. and Sierra, G., Nucl. Phys, 
{\bf 319} (1989) 155; ibid. B {\bf 330 } (1990) 347;
Phys. Lett., {\bf 220} (1989) 142; 
 C\'omez, C. and Sierra, G., Nucl. Phys., {\bf 352} (1991) 791.

\item[{[9]}]  Frenkel, B, and Reshetikhin, N, Yu., 
   Commun. Math, Phys. {\bf 146} (1992) 1.

\item[{[10]}]  Turaev, V. G., 
Inven. Math., 92 (1988), 527.

\item[{[11]}]  Akutsu, Y. and Wadati, M., 
J. Phys. Soc. Japan. {\bf 56} (1987) 839-842.

\item[{[12]}] Bazhanov, V. V. and Stroganov .Y. G., 
Theoret. Math. Fiz.,  
{\bf 62}(1985), 253-260.  

\item[{[13]}] M.L. Ge and K. Xue, 
J. of Phys. A: Math. \& Gen., {\bf 26} (1993), 281.  

\item[{[14]}] Gustav W. Delius, Mark D. Gould, Yao-Zhong Zhang, 
Nucl. Phys., {\bf B432} (1994) 377.

\item[{[15]}] Anthony J. Bracken, Mark D. Gould, Yao-Zhong Zhang and 
Gustav W. Delius, 
J. Phys., {\bf A27 }(1994), 6551.

\item[{[16]}] Sun, X. D., Wang, S. K. and Wu, K., 
Six-vertex type solutions of the colored Yang-Baxter equation, to  
appear in J. of Math. Phys..

\item[{[17]}] Fan, C. and  Wu, F. Y., 
Phys. Rev. {\bf B2}(1970), 723.

\item[{[18]}] Murakami, J., 
 A state model for the multi-variable Alexander polynomial 1990, 
preprint, Osaka University;

\item[{[19]}]  Cuerno,R., G\'omez,C., L\'opez, E.  and  Sierra,G., 
Phys. Lett. B., {\bf 307} (1993), 56-60.

\item[{[20]}] Murakami, J., 
Int. Jour. Mod. Phys., {\bf A7} Suppl. {\bf 1b} (1992) 765.

\item[{[21]}] Ruiz-Altaba, M, Phys. Lett., {\bf 277} (1992) 326.

\item[{[22]}] Wu Wen-tsun, Scientia Sinica, {\bf 21} (1978), 157-179.

\end{enumerate}
\vskip1cm

{\bf Shi-kuh Wang,  Institute of Applied Mathematics,
Chinese Academy of Sciences, Beijing 100080, P. R. of China.}

{\bf\sl  E-mail address: xyswsk@sunrise.pku.edu.cn}

\end{small}
\end{document}